\documentclass[10pt]{article}
\usepackage[utf8]{inputenc}
\usepackage[english]{babel}
\usepackage{amsmath,amsthm,amsfonts}
\usepackage{graphicx}
\usepackage{bm}
\usepackage{tikz}
\usepackage{multirow}
\usepackage{color}

\graphicspath{ {./figs/} }

\renewcommand{\div}{\mathop{\rm div}\nolimits}

\title{Upscaling method for problems in perforated domains with non-homogeneous boundary conditions on perforations using Non-Local Multi-Continuum method (NLMC)}

\author{
Maria Vasilyeva \thanks{Institute for Scientific Computation, Texas A\&M University, College Station, TX 77843-3368 \& Department of Computational Technologies, North-Eastern Federal University, Yakutsk, Republic of Sakha (Yakutia), Russia, 677980. Email: {\tt vasilyevadotmdotv@gmail.com}.}
\and
Eric T. Chung \thanks{Department of Mathematics,
The Chinese University of Hong Kong (CUHK), Hong Kong SAR. Email: {\tt tschung@math.cuhk.edu.hk}.}
\and
Wing Tat Leung
\thanks{Department of Mathematics, Texas A\&M University, College Station, TX 77843-3368, USA.}
\and
 Yating Wang
 \thanks{Department of Mathematics, Texas A\&M University, College Station, TX 77843-3368, USA.}
 \and
Denis Spiridonov
\thanks{ North-Eastern Federal University, Yakutsk, Republic of Sakha (Yakutia), Russia, 677980.}
\and
}

\begin{document}

\maketitle

\begin{abstract}
In this paper, we present an upscaling method for problems in perforated domains with non-homogeneous boundary conditions on perforations. Our methodology is based on the recently developed Non-local multicontinuum method (NLMC). The main ingredient of the method is the construction of suitable local basis functions with the capability of capturing multiscale features and non-local effects. We will construct multiscale basis functions for the coarse regions and additional multiscale basis functions for perforations, with the aim of handling non-homogeneous boundary conditions on perforations. We start with describing our method for the Laplace equation, and then extending the framework for the elasticity problem and parabolic equations. The resulting upscaled model has minimal size and the solution has physical meaning on the coarse grid. We will present numerical results (1) for steady and unsteady problems, (2) for Laplace and Elastic operators, and (3) for Neumann and Robin non-homogeneous boundary conditions on perforations. Numerical results show that the proposed method can provide good accuracy and provide significant reduction on the degrees of freedoms. 
\end{abstract}

\section{Introduction}

%Main points:
%\begin{itemize}
%\item We consider Laplace, elastic and parabolic problems in  perforated domain with non-homogeneous boundary conditions on perforations.
%\item We construct upscaled model using nlmc method with additional basis for handling  non-homogeneous boundary conditions on perforations.
%\item We present results (1) for steady and unsteady problems, (2) for Laplace and Elastic operators, and (3) for Neumann and Robyn non-homogeneous boundary conditions on perforations.
%\end{itemize}

In this paper, we will develop a multiscale method for solutions of problems in perforated domains without scale separation and in the presence of non-homogeneous boundary conditions on perforations. 
Mathematical modelling for the problems in perforated domains is important in many real-world applications. 
These applications include fluid flow in porous media, mechanical processes in composite materials, and so on.
Non-homogeneous boundary conditions on perforations can occur for reactive flow through porous media. These problems have great importance for a lot of applications in physics, chemistry, geology, and biology \cite{hornung1991diffusion}.
For the problems described as idealized periodic domains, two-scale homogenization method can be used for the construction of macroscale models  \cite{allaire1992homogenization, sanchez1980non, bakhvalov1984homogenization}. 
The homogenization techniques with the presence of a chemical reaction on perforations (solid grain interface for a porous media) lead to the additional reaction term in the macroscale problem \cite{allaire2010homogenization, korneev2016sequential, battiato2011hybrid, battiato2011hybrid}. 

Solutions of the problems in perforated domains have multiscale nature. Direct numerical solutions can lead to very large systems since a sufficiently fine computational mesh is needed to resolve the irregular boundaries of perforations as well as oscillations in the solutions. Thus, there are needs for some more efficient algorithms or multiscale methods. 
There are in literature a variety of multiscale approaches including 
the Heterogeneous Multiscale Method (HMM), the Mulitiscale Finite Element Method (MsFEM) and Mulitiscale Finite Volume Method (MsFVM) \cite{Henning09, le2014msfem, CR13, eh09, hw97, weinan2007heterogeneous, jenny2005adaptive, lunati2006multiscale}. 
In our previous works in solving problems in perforated domains, we use the multiscale model reduction technique based on the Generalized Mulitiscale Finite Element Method (GMsFEM) \cite{CELV2015, chung2017online, chung2016mixed, chung2018mstransport, chung2017conservativestokes}.
The GMsFEM is a general multiscale procedure, in which the model reduction is based on some local multiscale basis functions. The basis functions are constructed using local spectral decomposition \cite{egh12, EGG_MultiscaleMOR, chung2015residual, chung2016adaptive, eglp13oversampling, randomized2014}. 
The main idea behind the construction of the multiscale space is to design appropriate snapshot spaces and determine an appropriate
local spectral problem to select important modes in the snapshot space.
GMsFEM has been designed for many applications, for example, elasticity, thermoelasticity, poroelastic problems, wave propagation and so on \cite{ElasticGMsFEM, vasilyeva2016tp, MixedGMsFEM, chung2015wavefrac, brown2016generalized}.

In this paper, our goal is to develop a coarse-grid upscaled model, where the coarse grids do not have to align with perforations and the perforated domains do not need to have scale separation.  
For the construction of our upscaling method, we will use the recently developed Non-local multicontinuum method (NLMC) \cite{chung2017non}. In the work \cite{chung2017non}, we developed an upscaling method for flow problems in fractured porous media.
Upscaled model is directly related to the well-known multi-continuum approaches, which have been commonly used in approximating subgrid effects for flow and transport in fractured media \cite{chung2017coupling, warren1963behavior, barenblatt1960basic, douglas1990dual, akkutlu2015multiscale, akkutlu2017multiscale}.
Our Non-local multicontinuum upscaling \cite{chung2017non} is based on the recently developed  Constraint Energy Minimizing Generalized Multiscale Finite Element Method (CEM-GMsFEM) \cite{chung2017constraint}. In CEM-GMsFEM, one constructs multiscale basis functions so that they can capture long channelized effects and at the same time localizable. 
The construction of the multiscale space is based on an auxiliary space, which consists of eigenfunctions of local spectral problems. These auxiliary functions correspond to small (contrast-dependent) eigenvalues and represent the channels (high-contrast features). Using the auxiliary space, a constraint energy minimization problem is used to construct the required multiscale spaces. Due to a localization property, the minimization is performed locally in an oversampling domain, which is a few coarse elements larger than the target coarse block. 
Using the multiscale basis functions, a non-local upscaled model is then constructed. 

%The constraints allow handling non-decaying components of the local minimizers.

In this paper, we extend the NLMC upscaling approach for problems in perforated domains. We will construct an upscaled model using NLMC method with additional basis for handling  non-homogeneous boundary conditions on perforations. We will consider the Laplace, elastic and parabolic problems and show that our presented upscaling method can be applied to all these problems.
We consider several numerical examples: (1) steady and unsteady problems, (2) Laplace and Elastic operators, and (3) Neumann and Robin non-homogeneous boundary conditions on perforations. We show that one can achieve a good accuracy with a very few degrees of freedom.  %We discuss it in our numerical results.
In addition, the upscaled solution has a physical meaning in the coarse grid level. 

The paper is organized as follows. In Section~\ref{sec:problem}, we present the problems under consideration and their fine-scale approximations.
In Section~\ref{sec:method}, we give the constructions of our upscaling method. Numerical results are shown in Section~\ref{sec:num1} and Section~\ref{sec:num2}.
Finally, a conclusion is given in Section~\ref{sec:con}.

%We start with describing upscaling method for Laplace and elasticity problems in Section 1 and 2. In Section 2, we describe method for Laplace and elasticity problems, then extend framework for solution of the parabolic problem and show a relation between multi-continuum approaches and classic coarse grid models for reactive diffusive flow.  
%We present numerical results in Section 3 for model steady state model problems in perforated domain with non-homogeneous Neumann boundary conditions on perforations. 
%In Section 4, we perform numerical investigation of the merhod for time-dependent problem with non-homogeneous Neumann and Robyn boundary conditions on perforations. In the end, we present Conclusions.

\section{Problem formulation and fine grid approximation}\label{sec:problem}

\begin{figure}[!h]
\centering
\includegraphics[width=0.5 \textwidth]{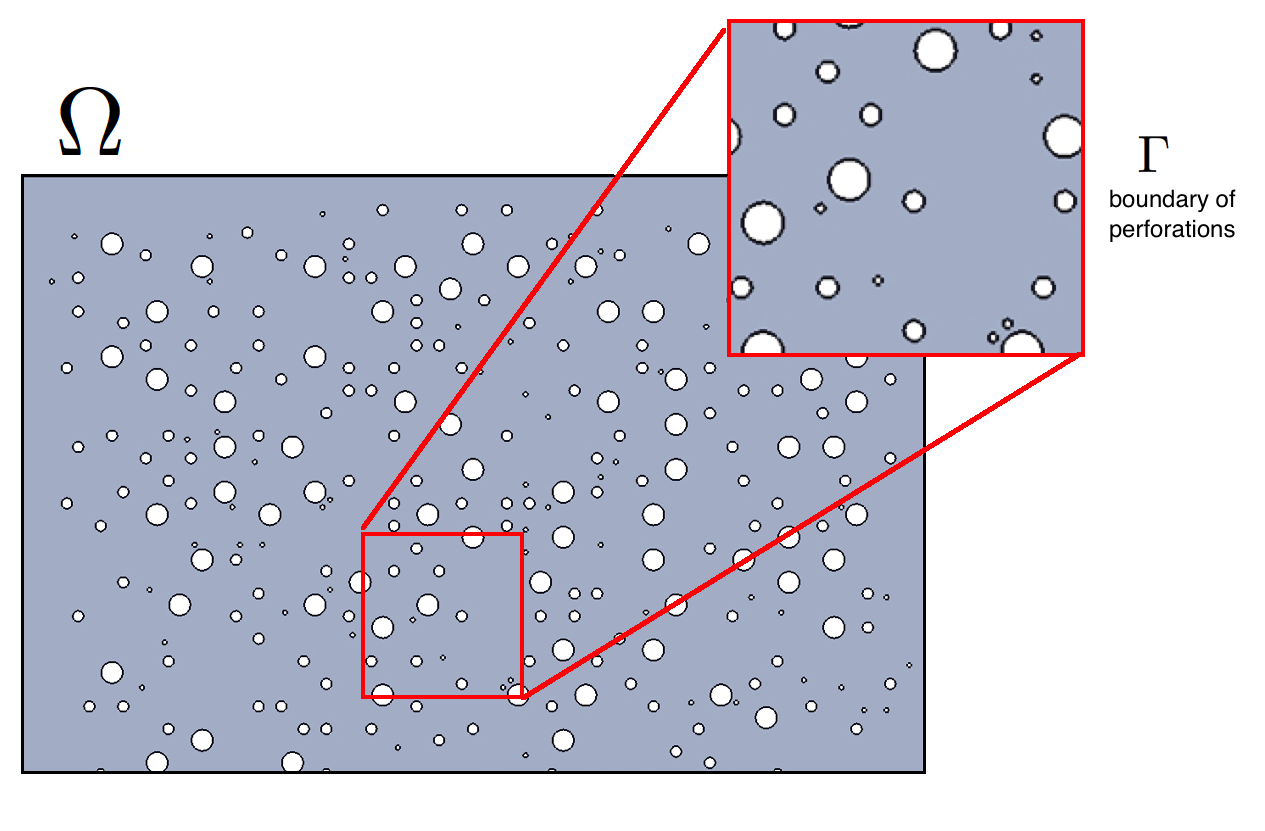}
\caption{An illustration of a perforated domain.}
\label{fig:domain}
\end{figure}

In this section, we will present the mathematical models under consideration,
and their standard fine scale approximations. Let $\Omega$ be a perforated domain (see Figure \ref{fig:domain} for an illustration).
% laplace model
We use $\Gamma$ to denote the boundary of the perforations,
and define $\partial \Omega \backslash \Gamma = \Gamma_D \cup \Gamma_N$ as the other part of the boundary of $\Omega$. 
We start our presentation with the Laplace problem in $\Omega$:
\begin{equation}
\label{eq:mm}
-\nabla \cdot (k \nabla u) = f, \quad \in \Omega,
\end{equation}
with a non-homogeneous Neumann boundary condition on the boundary of perforations $\Gamma$:
\begin{equation}
\label{eq:pb}
-k \nabla u \cdot n = g, \quad x \in \Gamma,
\end{equation}
and the following boundary conditions on $\Gamma_D \cup \Gamma_N$:
\begin{equation}
\label{eq:pb1}
u = 0, \quad x \in \Gamma_D,
\quad \text{ and } \quad
-k \nabla u \cdot n = 0, \quad x \in \Gamma_N.
\end{equation}
Here, $n$ denotes generically a unit normal vector for $\Gamma$ and $\Gamma_N$, $f$ denotes a given source
and $k$ is a heterogeneous coefficient. 

% elastic model
The second problem under consideration is the following elasticity problem in perforated domain $\Omega$:
\begin{equation}
\label{eq:mm2}
-\nabla \cdot \sigma(u) = f, \quad \in \Omega,
\end{equation}
where $u=(u_x,u_y)$, 
\[
\sigma (u) = \lambda \div u \, \mathcal{I} + 2 \mu \varepsilon (u),
\quad
\varepsilon (u) = \frac{1}{2} (\nabla u + (\nabla u)^T).
\]
Here $\varepsilon$ and $\sigma$ are the strain and stress tensors, $f$ is a given source vector, $\lambda$ and $\mu$ are the Lam\'{e} parameters.
We impose the above problem with a non-homogeneous Neumann boundary conditions on $\Gamma$:
\begin{equation}
\label{eq:pb2}
-\sigma \, n = g, \quad x \in \Gamma,
\end{equation}
and the following boundary conditions on $\Gamma_D \cup \Gamma_N$:
\[
u = 0, \quad x \in \Gamma_D,
\quad \text{ and } \quad
-\sigma \, n = 0, \quad x \in \Gamma_N.
\]

% approximation
For the numerical solution, we need a fine grid that resolves all perforations. We will use the standard piecewise linear
finite element space $V_h$ and have following variational formulation: find $u_h \in V_h$ such that
\begin{equation}
\label{eq:vf}
a(u_h, v) = l(v), \quad \forall v \in V_h,
\end{equation}
where  $l(v) = \int_{\Omega} f \, v\, dx + \int_{\Gamma} g \, v\, ds,$ and 
\begin{itemize}
\item for Laplace problem: 
$
a(u, v) = \int_{\Omega} k \nabla u \cdot \nabla v \, dx, 
$
\item for elasticity problem: 
$
a(u, v) = \int_{\Omega} \sigma(u) : \varepsilon(v) \, dx.
$
\end{itemize}
Assume that $V_h = \text{span} \{ \phi_i \}$, where $\phi_i$ is the standard linear nodal basis function, and  $u_h = \sum_i u_i \phi_i$.
Then we can write the above variational problem in matrix form
\begin{equation}
\label{eq:mf}
A U = b,
\end{equation}
where $U = (u_j)$, $b_j = \int_{\Omega} f \, \phi_j \, dx + \int_{\Gamma} g \, \phi_j \, ds$, and 
\begin{itemize}
\item for Laplace problem: 
$
A = \{a_{ij}\}, \quad a_{ij} = \int_{\Omega} k \nabla \phi_i \cdot \nabla \phi_j \, dx,
$
\item for elasticity problem: 
$
A = \{a_{ij}\}, \quad a_{ij} = \int_{\Omega} \sigma(\phi_i) : \nabla \varepsilon(\phi_j) \, dx.
$
\end{itemize}
Therefore, on fine grid,  we have a linear system of equations with size $DOF_f = N_f$ for Laplace problem, and  $DOF_f = 2 N_f$ for elasticity, where $N_f$ is the number of vertices
in the interior of $\Omega$ and on $\Gamma_N \cup \Gamma$.

Finally, we will consider the following time-dependent problem
\begin{equation}
\label{t-nlmc1}
c  \frac{\partial u}{\partial t} - \nabla \cdot (k \nabla u)  = f,
\end{equation}
subject with the boundary conditions (\ref{eq:pb}) and (\ref{eq:pb1}) as well as a suitable initial condition. 
%with non-homogeneous Neumann boundary conditions on perforations boundary  ($-k \nabla u \cdot n  = g$).
The fine grid approximation can be written as
\[
S \frac{u^{n+1} - u^n}{\tau} + A u^{n+1}  = b,
\]
where we use a stable implicit scheme for approximation in time, $S$ is the mass matrix, $\tau$ is the given time step and
\[
S = \{s_{ij}\}, \quad s_{ij} = \int_{\Omega} c \phi_i \phi_j \, dx, \quad 
A = \{a_{ij}\}, \quad a_{ij} = \int_{\Omega} k \nabla \phi_i \cdot \nabla \phi_j \, dx, \quad
b_j = \int_{\Omega} q \, \phi_j \, dx + \int_{\Gamma} g \, \phi_j \, ds.
\]
For this problem, one also needs to solve a system of size $DOF_f = N_f$.

\section{Coarse grid upscaling using NLMC method}\label{sec:method}

The main idea of the NLMC method is to compute the upscaled multi-continuum coefficients by some appropriate multiscale basis functions. In the work \cite{chung2017non}, where flows in fractured porous media are considered, each multi-continuum coefficient corresponds to a fracture network and is able to capture both local and non-local effects. The construction of the basis functions is based on the CEM-GMsFEM \cite{chung2017constraint}. In particular, basis functions are solutions of an energy minimization problem subject to some appropriate orthogonality conditions obtained by local spectral problems. Moreover, a localization property is proved, and this leads to the local computations of basis functions on oversampled regions obtained by enlarging the target coarse element by a few coarse grid layers, which depend weakly on the contrast. Besides, convergence theory is established which states that the method is convergent with respect to coarse mesh size and independent of the medium properties. For more details of NLMC method applied to flows in fractured media, see \cite{chung2017non}. We remark that the NLMC method is related to the dual porosity (multi-continuum) approach.

%In this method, we follow a general concept of spectral basis functions in GMsFEM and apply simplified approach, where we define constant functions in perforations and the background medium. For construction of the multiscale basis functions, we solve local problems in an oversampled local domains subject to the constraint that the local solution vanishes in other continua except the one for which it is formulated for.

In this paper, we present the construction of an accurate coarse grid model for problems in perforated domains using NLMC method. We will use a new and simplified approach for the construction of local basis functions without the need of local spectral problems. For each local region, we will construct a basis function corresponding to the background medium and some basis functions for each the perforations in the target local region. We will solve local energy minimization problems on appropriate oversampled regions subject to the constraint that the solution vanishes in other continua except the one for which it is formulated. 
%We construct simplified basis for perforated domain which has spatial decay property in local domains and separate background medium and perforations. This technique is similar to the dual porosity (multi-continuum) approach that used for flow problems in fractured porous media, where we have coupled system of equations for flow in porous matrix and fracture networks.
We will develop two types of basis functions. For the first type, we will use one single basis function for all perforations within a target coarse element. 
For the second type, we construct a basis function for each perforation within a target coarse element. We remark that while the second type can be more accurate, the first type is more economical as the number of perforations can be very large and can lead to a significant increase in the coarse grid model size.

%constant withing each separate perforation in each coarse cell and a constant for background medium. For second type, we use one additional basis for all perforations in each coarse cell, that is better for problems in perforated domains because number of perforations can be very large and can lead to the significant increasing of the coarse grid model size.

\begin{figure}[h!]
\centering
\includegraphics[width=0.7 \textwidth]{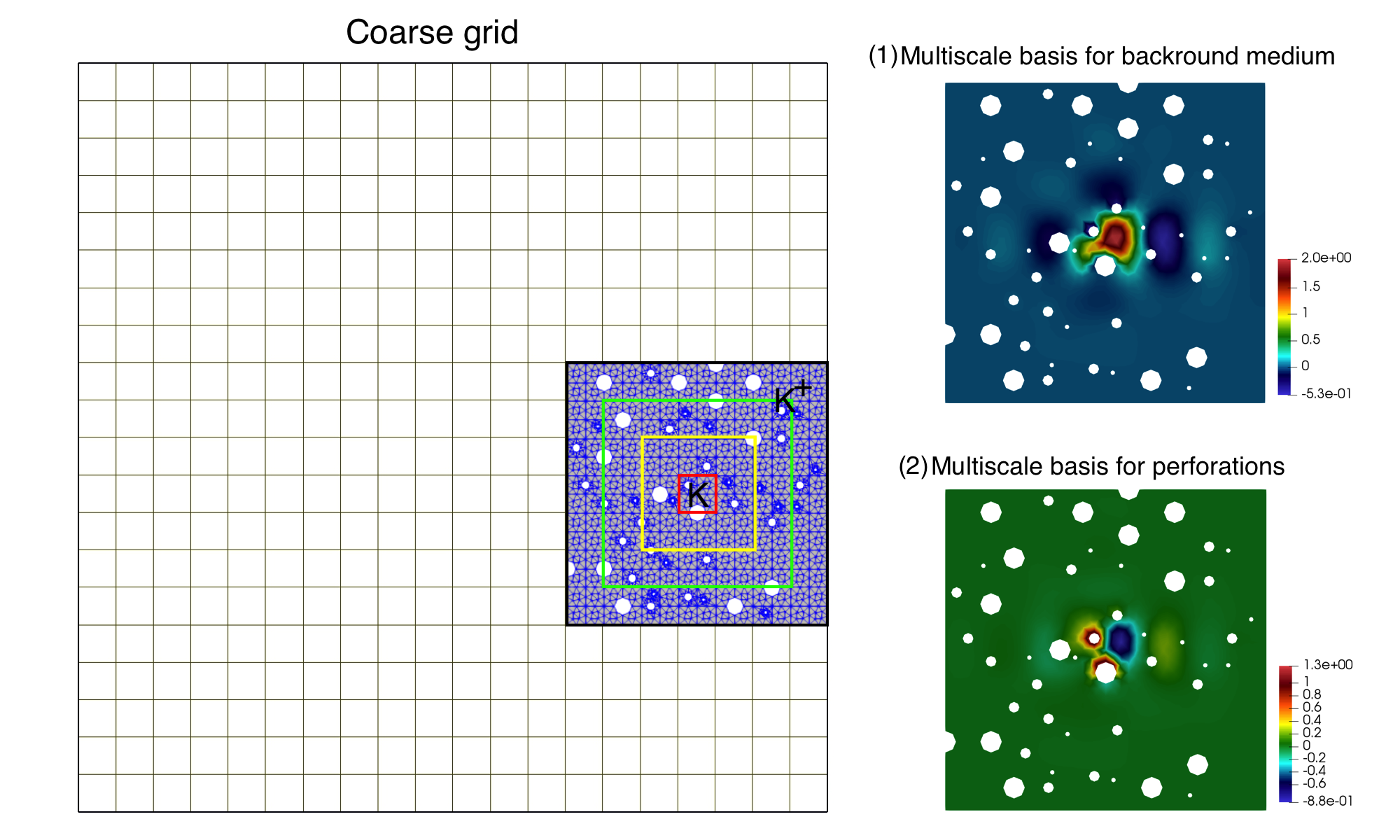}
\caption{Left: A coarse block $K$ and its oversampled region $K^3$.
Right: Multiscale basis functions for the coarse block $K$ with support on $K^3$.}
\label{fig:msbf20}
\end{figure}

We will now give the details of the constructions. 
We consider a coarse partition $\mathcal{T}^H$ of the domain $\Omega$.
Let $K_i \in\mathcal{T}^H$ be the $i$-th coarse block and 
let $K^+_i$ be the corresponding oversampled region obtained by expanding the coarse block $K_i$ by several coarse grid layers
(see Figure \ref{fig:msbf20} for an illustration).
For the ease of presentation, we use $K^n$ to denote the oversampled region obtained by expanding $K$ by $n$ coarse grid layers. 
In Figure \ref{fig:msbf20}, we give an illustration of the basis functions computed for a coarse block $K$, and these basis functions
correspond to the background medium and perforations,
and they have supports on $K^3$.

% Laplace basis

We let $N_c$ be the number of elements in $\mathcal{T}^H$. 
For a given region $S \subset \Omega$ which is a connected union of coarse grid elements,
we define $V_h(S)$ as the space of functions in $V_h$ with support in $S$ and with zero trace on $\partial S \backslash (S \cap \Gamma)$. 
To define the two types of basis functions, we need two types of constraints. 
We remark that these constraints are constructed for the basis functions corresponding to the target coarse element $K_i$. 
We also remark that the following construction is for the Laplace operator. 

\textbf{Type 1 basis functions.}
We will define the Type 1 basis functions. These are functions in the space $V_h(K_i^+)$, where $K_i^+$ is an oversampled region for $K_i$.
For each coarse block $K_j \subset K_i^+$, we define $\gamma_j = K_j \cap \Gamma$, which is the set of boundaries of the perforations in the block $K_j$.
We will construct two basis functions $\psi^i_0$ and $\psi^i_1$, which minimize the energy $a(\psi,\psi)$ restricted in $V_h(K_i^+)$
and satisfy the constraints described below:

%The purpose is to define a constrained space $V_i^{con} \subset V_h(K_i^+)$, where $K_i^+$ is an oversampled region for $K_i$.
%For each coarse block $K_j \subset K_i^+$, we define $\gamma_j = K_j \cap \Gamma$, which is the set of boundaries of the perforations in the block $K_j$.
%The space $V_i^{con}$ is spanned by two functions $\phi^i_0$ and $\phi^i_1$, which correspond to the background medium and the perforations respectively. 
%These two functions are defined so that, for each $K_j \subset K_i^+$, the following conditions hold:

%Let $\gamma_j = K_j \cap \Gamma$ is the all perforations' boundary inside coarse cell $K_j \in K_i^+$.
%As a simplified auxiliary space, we use $V^{aux}(K_i) = \text{span} \{ \phi^i_{l} \}$, where $i = \overline{1,N_c}$ ($N_c$ is the number of coarse grid cells) and $l = 0,1$ ($L_i$ is the number of perforations in $K_i$). For functions $\phi^i_{l}$, we have following conditions $\forall K_j \in K_i^+$

\begin{itemize}
\item $\psi_0^i$ (backround medium) :
\[
\int_{K_j} \psi_0^i \, dx = \delta_{i,j}, \quad
\int_{\gamma_j} \psi_0^i \, ds = 0.
\]
\item $\psi_1^i$ (perforation):
\[
\int_{K_j} \psi_1^i \, dx = 0, \quad
\int_{\gamma_j} \psi_1^i \, ds = \delta_{i,j}.
\]
\end{itemize}

In the above definitions, $\delta_{i,j}$ is the Dirac delta function. 
We remark that the basis functions $\psi^i_0$ and $\psi^i_1$ correspond to the background medium and the perforations respectively. 
In the case that $K_i$ does not contain perforations, we only construct one basis function $\psi_0^i$.
Note that, the function $\psi_0^i$ has an average value one in the background medium and zero average on perforations, and the second function $\psi_1^i$ has an average zero in the background medium and average one on perforations.

\textbf{Type 2 basis functions.}
We will define the Type 2 basis functions. These are functions in the space $V_h(K_i^+)$, where $K_i^+$ is an oversampled region for $K_i$.
We first write $\Gamma = \cup_{l = 1}^L \Gamma^{(l)}$, where $\Gamma^{(l)}$ is the $l$-th perforation and $L$ is the total number of perforations.
For each coarse block $K_j \subset K_i^+$, we define $\gamma_j^{(l)} = K_j \cap \Gamma^{(l)}$, which is the boundary of the $l$-th perforation in the block $K_j$.
We let $L_i$ be the number of perforations in $K_i$. 
We will construct a set of basis functions $\psi^i_0$ and $\psi^i_k$ ($k=1,2,\cdots, L_i$), which minimize the energy $a(\psi,\psi)$ restricted in $V_h(K_i^+)$
and satisfy the constraints described below:

%\textbf{Auxiliary space of Type 2. }  
%Let $\Gamma = \cup_{l = 1}^L \Gamma^{(l)}$, where $\Gamma^{(l)}$ - is the one perforation and $L$ is the total number of perforations.
%$\gamma^{(l)}_j = K_j \cap \Gamma^{(l)}$ is the perforations inside coarse cell $K_j \in K_i^+$ and  $L_j$ is the number of perforations in $K_j$. 
%Then, as a simplified auxiliary space, we use  $V^{aux}(K^+_i) = \text{span} \{ \phi^i_{l} \}$, where $i = \overline{1,N_c}$ ($N_c$ is the number of coarse grid cells) and $l = \overline{0, L_i}$ ($L_i$ is the number of perforations in $K_i$). For functions $\phi^i_{l}$, we have following conditions $\forall K_j \in K_i^+$
\begin{itemize}
\item $\psi^i_{0}$ (background medium) :
\[
\int_{K_j} \psi_0^i \, dx = \delta_{i,j}, \quad 
\int_{\gamma^{(l)}_j} \psi_0^i \, ds = 0, \quad  l=\overline{1, L_j}.
\]
\item  $\psi_m^i$ ($m$-th perforation in $K_i$):
\[
\int_{K_j} \psi_m^i \, dx = 0, \quad 
\int_{\gamma^{(l)}_j} \psi_m^i \, ds = \delta_{i,j}\delta_{m,l}, \quad  l=\overline{1, L_j}.
\]
\end{itemize}
%This conditions, we use as constrains for multiscale basis functions construction.

%\textbf{Basis functions.}  
%For construction of the multiscale basis functions, we ensure constrain energy minimizing property
%\[
%\psi^i_l = \text{argmin} \{ a(\psi, \psi) \, | \, \psi \in V_0^{aux}(K^+_i) \},
%\]
%where $V_0^{aux}(K^+_i) = \{ v \in V^{aux}(K^+_i): v = 0 \text{ on } \partial K^+_i / \Gamma \}$.

We remark that, we can obtain the Type 1 basis functions $\psi^i_m$ ($m=0,1$) by solving the following local problem on $K_i^+$:
\begin{equation}
\label{eq:basis}
\begin{split}
&a(\psi^i_m, v)
+ \sum_{K_j \subset K_i^+} \lambda_0 \int_{K_j} v \, dx
+ \sum_{K_j \subset K_i^+} \lambda_1 \int_{\gamma_j} v \, ds = 0, \\
& \int_{K_j} \psi^i_m \, dx = \delta_{i,j}, \quad \forall K_j \in K_i^+, \\
& \int_{\gamma_j} \psi^i_m \, dx = \delta_{i,j} \delta_{m,l}, \quad \forall \gamma_j = K_j \cap \Gamma,
\end{split}
\end{equation}
with zero Neumann boundary conditions on perforations $\Gamma$ and zero Dirichlet boundary conditions on $\partial K^+_i \backslash \Gamma$. Note that, we use a Lagrange multiplier to enforce the constraints. See Figure \ref{fig:msbf20} for an illustration of the Type 1 basis functions.
Similarly, we can construct the Type 2 multiscale basis functions by the method of Lagrange multipliers.

%In Figure \ref{fig:msbf20}, we depict a multiscale basis functions for Laplace problem in perforated domain with oversampled region $K^+_i = K^3_i$ (three oversampling coarse cell layers) on coarse mesh $20 \times 20$.

%
\textbf{Coarse grid system.}
Finally, we obtain our multiscale space
\[
V_{ms} = \text{span} \{ \psi^i_m \}
\]
using our multiscale basis functions. 
We use these local multiscale basis functions to construct the upscaled equation. Following \cite{chung2017non},
the resulting coarse-grid equation can be written in a following discrete form 
\begin{equation}
\label{nlmc1}
\sum_{j, n}  T^{mn}_{ij} \bar{u}^n_j  = \bar{q}^m_i, \quad 
T^{mn}_{ij} =  a(\psi^m_i, \psi^n_j), \quad 
\bar{q}^m_i =  l(\psi^n_j),
\end{equation}
for  cell $K_i$ and continuum $m$, where $m = 0,1$ for Laplace problem, $m=0$ is related to the background medium and $m=1$ is related to the perforations. 
Construction of the $T^{mn}_{ij}$ can be done in the offline stage as precalculations.
We remark that the upscaled solution is denoted as $(\bar{u}^n_j)$.

The implementation of the method is discussed as follows:
\begin{enumerate}
\item Calculation of the multiscale basis functions $\psi^i_{l}$ ($l = 0, 1,...,L_i$) for background medium and perforations by solution of the local problems in $K^+_i$ for each $i = \overline{1, N_c}$.

\item Generation of the projection matrix
\[
R^T = \left[ \psi_0^1, \ldots, \psi_{L_1}^{1}, \ldots, \psi_{0}^{N_c}, \ldots, \psi_{L_{N_c}}^{N_c} \right],
\]
where we understood $\psi_l^i$ as a column vector using its representation in the fine grid. 
%$\phi^i_{l}$ are used to denote the nodal values of each basis function defined on the fine grid.

\item Construction of the coarse grid system
\[
\underbrace{R A R^T}_{T} \bar{u} = \underbrace{R b}_{\bar{q}},
\]
and solve $T \bar{u} = \bar{q}$.
\end{enumerate}

Note that, if needed, we can reconstruct the downscale solution, $u_{ms} = R^T \bar{u}$. Our coarse grid solutions have physical meaning, which is the average value of the solution on each coarse cell and boundary of perforations thanks to the construction of the multiscale basis functions.

\textbf{NLMC upscaling for elasticity problem.}
% elasticity
For the construction of the multiscale basis functions for the elasticity problem (\ref{eq:mm2}), we use a similar approach.
We will present the construction for Type 2 basis functions. In particular,
we will construct a set of basis functions $\psi^{X,i}_l :=  (\psi^{X,i}_{x, l}, \psi^{X,i}_{y, l})$ 
and $\psi^{Y,i}_l := (\psi^{Y,i}_{x, l}, \psi^{Y,i}_{y, l})$, ($l=0,1,\cdots, L_i$), which minimize the energy $a(\psi,\psi)$ restricted in $V_h(K_i^+)$
and satisfy the constraints described below for all $K_j \subset K_i^+$. To reduce the repetitions, we will show the constraints for $\psi^{X,i}_l$,
and the constraints for $\psi^{Y,i}_l$ are defined analogously. 

\begin{itemize}
\item $\psi^{X,i}_{0}$ (background medium) :
\begin{eqnarray*}
&& \int_{K_j} \psi_{x,0}^{X,i} \, dx = \delta_{i,j}, \quad 
\int_{\gamma^{(l)}_j} \psi_{x,0}^{X,i} \, ds = 0, \quad  l=\overline{1, L_j}, \\
&& \int_{K_j} \psi_{y,0}^{X,i} \, dx = 0, \quad 
\int_{\gamma^{(l)}_j} \psi_{y,0}^{X,i} \, ds = 0, \quad  l=\overline{1, L_j}.
\end{eqnarray*}
\item  $\psi_m^{X,i}$ ($m$-th perforation in $K_i$):
\begin{eqnarray*}
&& \int_{K_j} \psi_{x,m}^{X,i} \, dx = 0, \quad 
\int_{\gamma^{(l)}_j} \psi_{x,m}^{X,i} \, ds = \delta_{i,j}\delta_{m,l}, \quad  l=\overline{1, L_j}, \\
&& \int_{K_j} \psi_{y,m}^{X,i} \, dx = 0, \quad 
\int_{\gamma^{(l)}_j} \psi_{y,m}^{X,i} \, ds = 0, \quad  l=\overline{1, L_j}.
\end{eqnarray*}
\end{itemize}

\begin{figure}[h!]
\centering
\includegraphics[width=1 \textwidth]{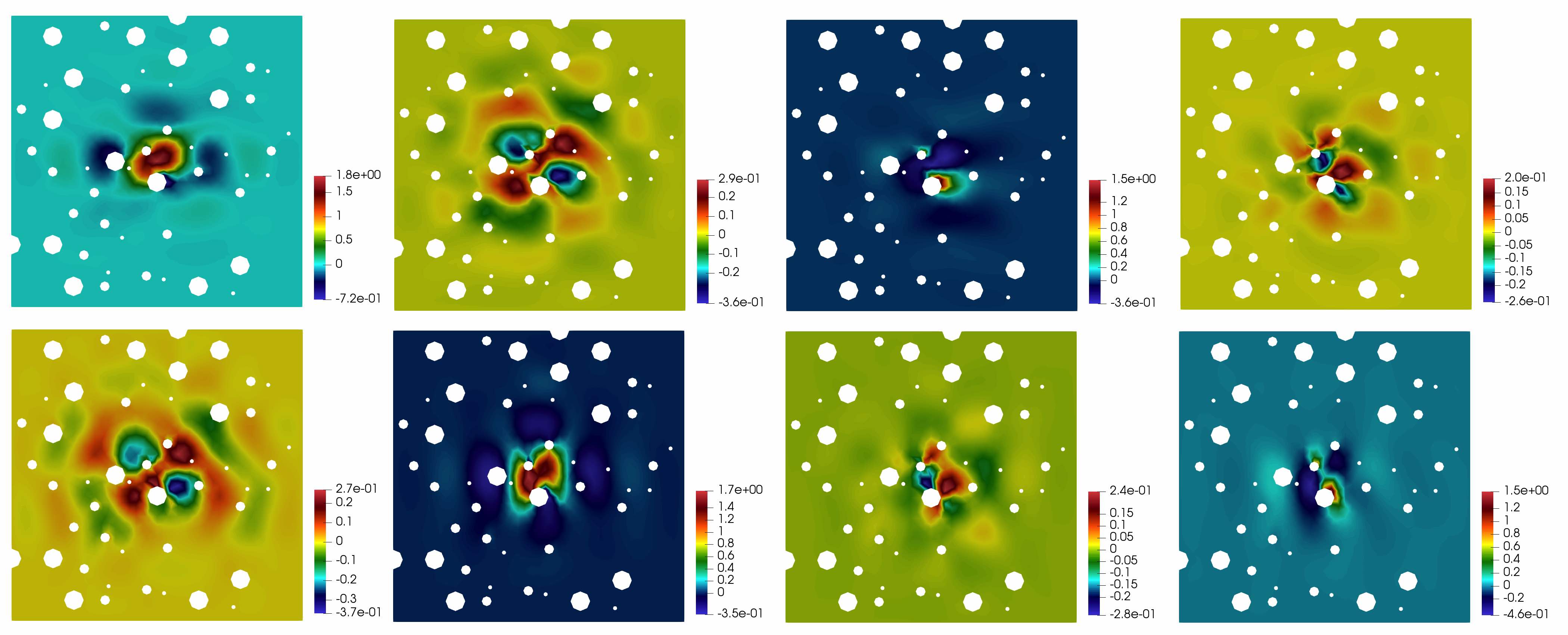}
\caption{Basis functions for elasticity problem computed on an oversampled region $K_i^+ := K_i^3$. 
First row: $\psi^{X,i}_l$. Second row: $\psi^{Y,i}_l$. Here, we have $l = 0,1,2,3$ (from left to right).}
\label{fig:msbf20xy}
\end{figure}

The construction of Type 1 multiscale basis functions are defined in a similar fashion. In Figure \ref{fig:msbf20xy},
we present a set of Type 2 basis functions $\psi^{X,i}_l$ and $\psi^{Y,i}_l$ with $l=0,1,2,3$
for a selected coarse block $K_i$ computed on $K_i^+ := K_i^3$.

% elastic basis

%For construction of the multiscale basis functions, we use similar local problem as before with elastic operator. Analogously, the multiscale space of Type 1 can be constructed (see Figure \ref{fig:msbf20xy} for multiscale basis functions for elasticity problem in perforated domain with oversampled region $K^+_i = K^3_i$ (three oversampling coarse cell layers) on coarse mesh $20 \times 20$).

Finally, we define the projection matrix for the elasticity problem as
\[
R = \binom{R_X}{R_Y}, \quad
R_X^T = \left[ \psi_0^{X,1}, \ldots, \psi_{L_1}^{X,1}, \ldots,  \psi_0^{X,N_c}, \ldots, \psi_{L_{N_c}}^{X,N_c} \right],
\quad
R_Y^T = \left[ \psi_0^{Y,1}, \ldots, \psi_{L_1}^{Y,1}, \ldots,  \psi_0^{Y,N_c}, \ldots, \psi_{L_{N_c}}^{Y,N_c} \right],
\]
where we understood the basis functions as column vectors using their representations in the fine grid. 
Then, we generate the coarse grid system as
$T \bar{u} = \bar{q}$ with
$T = R A R^T$, $\bar{q} = R b$ and $\bar{u} = (\bar{u}_X, \bar{u}_Y)$.

\textbf{Time-dependent problem.} 
Similar to the two problems discussed above, we can formulate the upscaled model for the time-dependent equation (\ref{t-nlmc1}) as follows:
\begin{equation}
\label{t-nlmc2}
M \frac{\bar{u}^{n+1} - \bar{u}^n}{\tau} + T \bar{u}^{n+1} = \bar{q},
\end{equation}
where the mass and stiffness matrices ($M = \{s^{mn}_{ij}\}$ and $T = \{a^{mn}_{ij}\}$) are defined as
\[
s^{mn}_{ij} =
\left\{\begin{matrix}
\sum_{j,n} s^{mn}_{ij}, & i = j, m = n, \\
0, & else
\end{matrix}\right. ,
\quad
a^{mn}_{ij} =
\left\{\begin{matrix}
- \sum_{j \neq i, n \neq m} a^{mn}_{ij}, & i = j, m = n\\
a^{mn}_{ij}, & else
\end{matrix}\right. ,
\]
with $s^{mn}_{ij} = s(\psi^m_i, \psi^n_j)$ and $a^{mn}_{ij} = a(\psi^m_i, \psi^n_j)$, and the basis functions $\{ \psi^m_i\}$ are defined in above for the Laplace problem. 
Note that, mass matrix is diagonal and by the properties of the constructed multiscale basis functions, we can directly calculate the mass matrix elements on the coarse grid as
$s^{mm}_{ii}= c_m |V_i| / \tau$ and right-hand side vector $\bar{q}^m_i = f_m |V_i| + h_m |V_i|$ ($h_0$ = 0 and $h_m = g$ for $m \neq 0$), where $|V_{0}| = |K_i|$ and $|V_{m}| = |\gamma_m|$ for $m \neq 0$.
Therefore, we can write our upscaled system as 
\begin{equation}
\label{t-nlmc3}
\left(
\frac{1}{\tau}
\begin{pmatrix}
M_b & 0 \\
0 & M_p
\end{pmatrix}
+
\begin{pmatrix}
T_b & T_{bp} \\
T_{pb} & T_p
\end{pmatrix}
\right)
\binom{\bar{u}_b^{n+1}}{\bar{u}_p^{n+1}}
=
\binom{\bar{q}_b}{\bar{q}_p}
+
\frac{1}{\tau}
\binom{M_b \bar{u}_b^n}{M_p \bar{u}_p^n},
\end{equation}
where $\bar{u} = (\bar{u}_b, \bar{u}_p)$, $\bar{u}_b$ and $\bar{u}_p$ are the average cell solution on coarse grid for background medium and for perforations.
Mass matrix is diagonal, stiffness matrix is non-local and provides a good approximation due to the coupled construction.

% Robyn
Finally, we consider a coarse-grid upscaled model for problem \eqref{t-nlmc1} with  Robin boundary conditions on perforations
 \[
-k \nabla u \cdot n  = \alpha (u - g), \quad x \in \Gamma.
\]
We can write the fine grid approximation as
\[
S \frac{u^{n+1} - u^n}{\tau} + A u^{n+1} + B u^{n+1}  = b,
\]
where $B$ is the boundary mass matrix 
\[
B = \{r_{ij}\}, \quad r_{ij} = \int_{\Gamma} \alpha \psi_i \psi_j \, ds, \quad 
b_j = \int_{\Omega} q \, \psi_j \, dx + \int_{\Gamma} \alpha g \, \psi_j \, ds.
\]
Similar as before, we have the following upscaled model
\[
M \frac{\bar{u}^{n+1} - \bar{u}^n}{\tau} + T \bar{u}^{n+1} + C  \bar{u}^{n+1}  = \bar{q},
\]
where $C = \{r^{mn}_{ij}\}$,
\[
r^{mn}_{ij} = 
\left\{\begin{matrix}
\sum_{j,n} r^{mn}_{ij}, & i = j, m = n, \\ 
0, & else,
\end{matrix}\right.
\]
with $r^{mn}_{ij} = r(\psi^m_i, \psi^n_j)$. This boundary mass matrix is diagonal and can be approximated by $r^{mm}_{ii}= \alpha h_m |V_i| $, where $h_0 = 0$ and for $m \neq 0$, we set $h_m = 1$ with $|V_{m}| = |\gamma_m|$. We note that, the summation of each row for the matrix $T$ is zero which ensures the mass conservation.

\section{Numerical results}

In this section, we will show numerical results to demonstrate the performance of our upscale method. We will consider the Laplace and the elasticity problems in
Section~\ref{sec:num1} and the time dependent problem in Section~\ref{sec:num2}.

\subsection{Numerical results for Laplace and elasticity problems}\label{sec:num1}

\begin{figure}[h!]
\centering
\includegraphics[width=0.32 \textwidth]{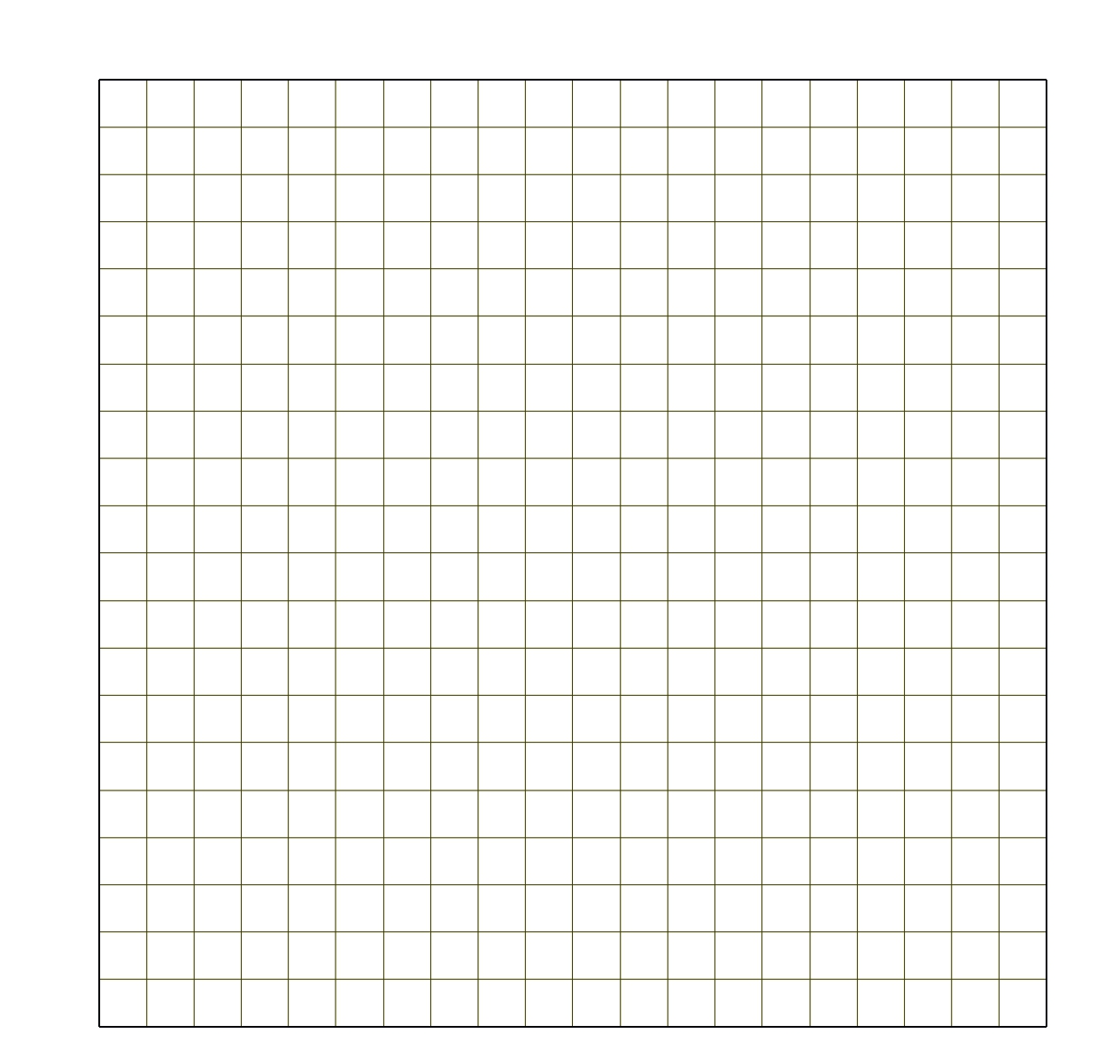}
\includegraphics[width=0.32 \textwidth]{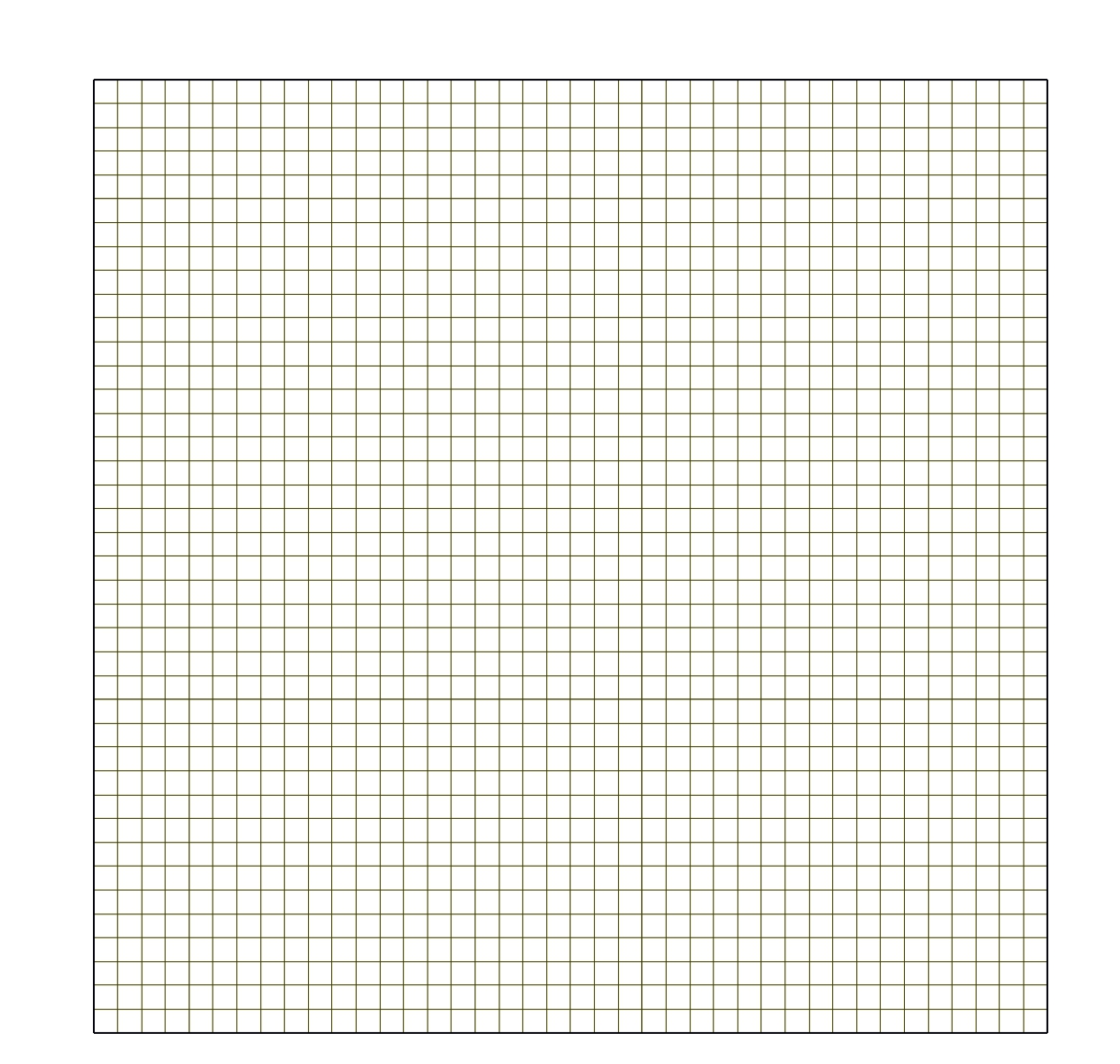}
\includegraphics[width=0.32 \textwidth]{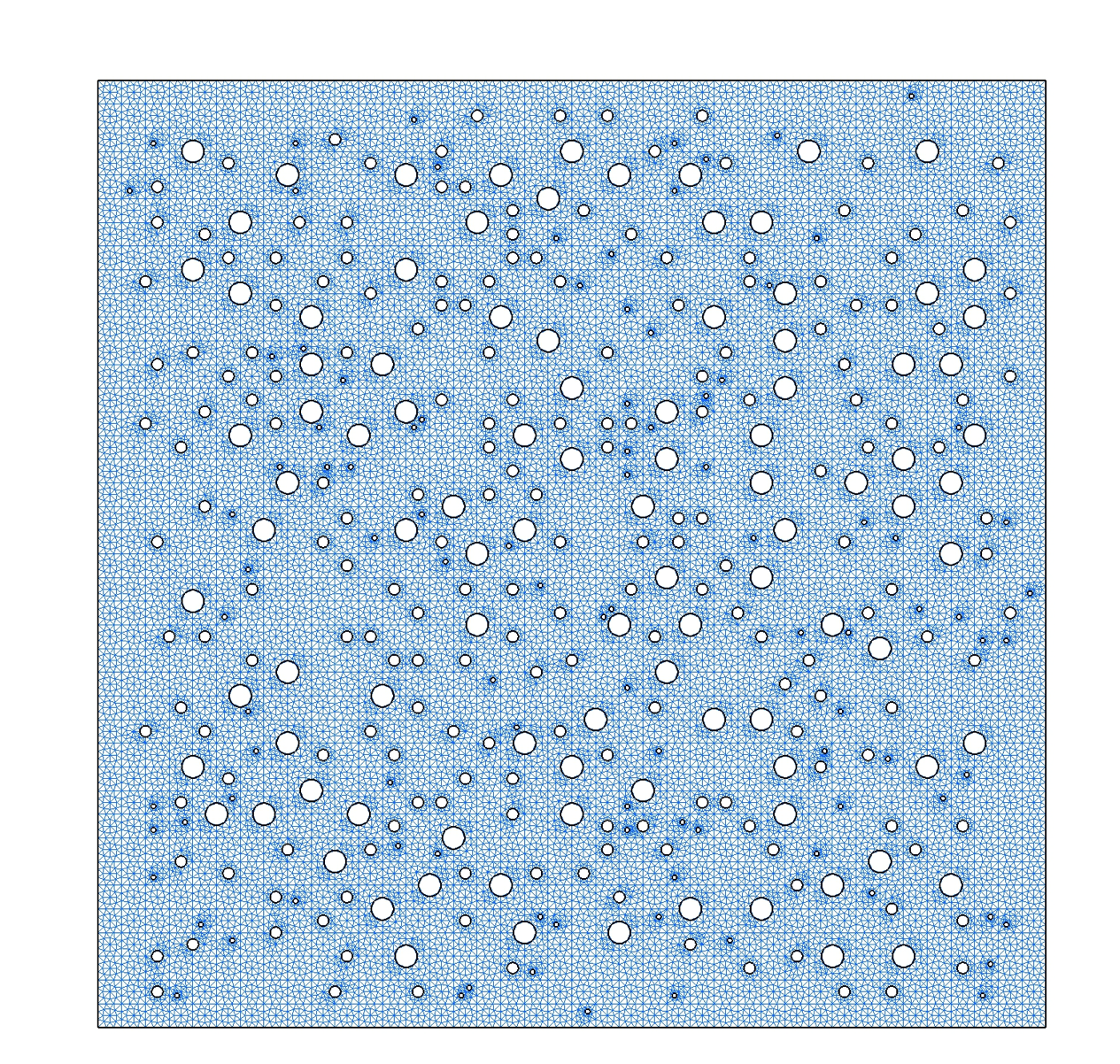}
\caption{Computational grids.  Left: Coarse grids with 400 ($20 \times 20$) cells. Middle:  Coarse grids with 1600 cells ($40 \times 40$). Right: Fine mesh with 15389 vertices. }
\label{fig:mesh}
\end{figure}

We present numerical results for Laplace and elasticity problems in domain $\Omega = [0, 1] \times [0, 1]$ with 400 perforations. These perforations are resolved on the fine grid using triangular cells. The coarse grid is uniform with rectangular cells.
In Figure \ref{fig:mesh}, we show computational coarse and fine grids. For numerical simulation, we consider two grids: (1) coarse grid $20 \times 20$ (400 cells) and (2) coarse grid $40 \times 40$ (1600 cells).
We use $DOF_c$ to denote problem size of the upscaled model. 

To compare the results, we use the relative $L^2$ error between coarse cell average of the fine-scale solution $\bar{u}_f$ and upscaled coarse grid  solutions $\bar{u}$
\begin{equation}
e_{L_2} = ||\bar{u}_f - \bar{u} ||_{L^2}, \quad
|| \bar{u}_f - \bar{u} ||^2_{L^2} =  
\frac{ \sum_K (\bar{u}^K_f - \bar{u}^K)^2 \, dx}{\sum_K (\bar{u}^K_f)^2 \, dx}, \quad 
\bar{u}^K_f = \frac{1}{|K|} \int_K u_f \, dx.
\end{equation}

% steady
\begin{figure}[h!]
\centering
\includegraphics[width=0.75 \textwidth]{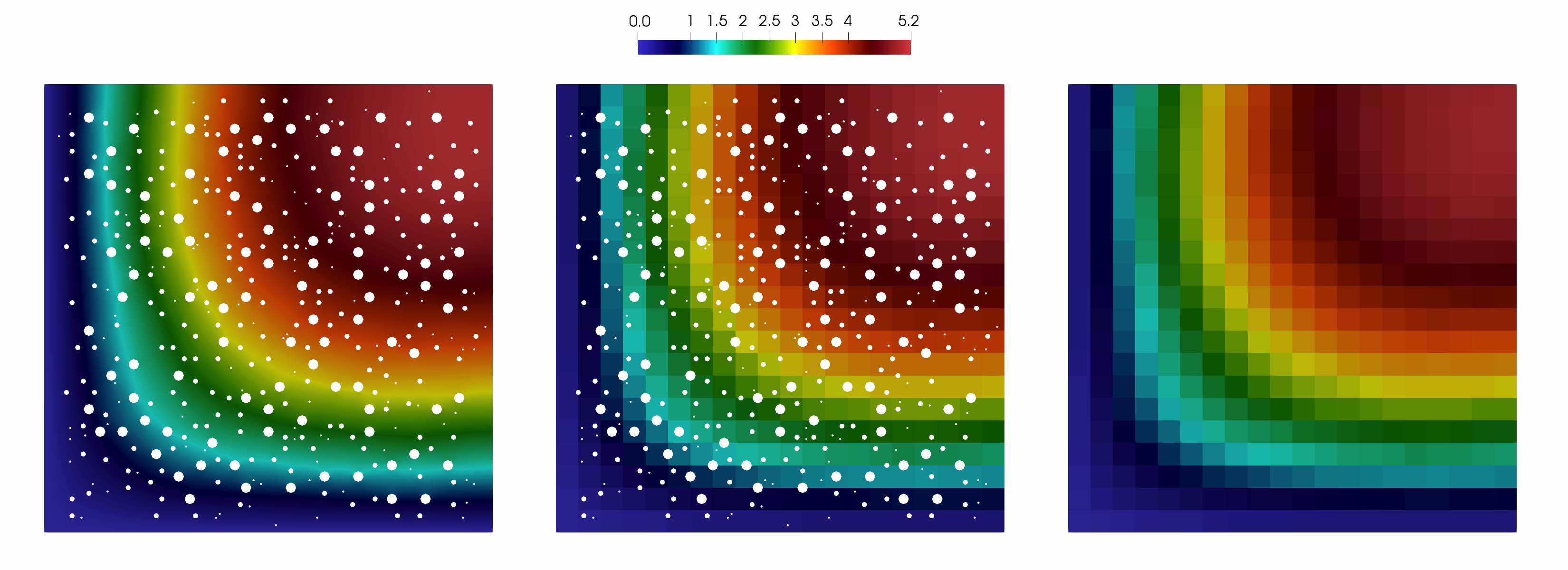}
\caption{Laplace problem in the perforated domain shown in Figure~\ref{fig:mesh}. 
Left: fine-scale solution using $DOF_f = 15389$. Middle: mean value of fine-scale solution on a $20 \times 20$ coarse mesh. Right: multiscale solutions computed on a $20 \times 20$ coarse mesh with $DOF_c = 961$ using $4$ oversampling layers in basis constructions.}
\label{fig:u20}
\end{figure}

% elastic
\begin{figure}[h!]
\centering
\includegraphics[width=0.75 \textwidth]{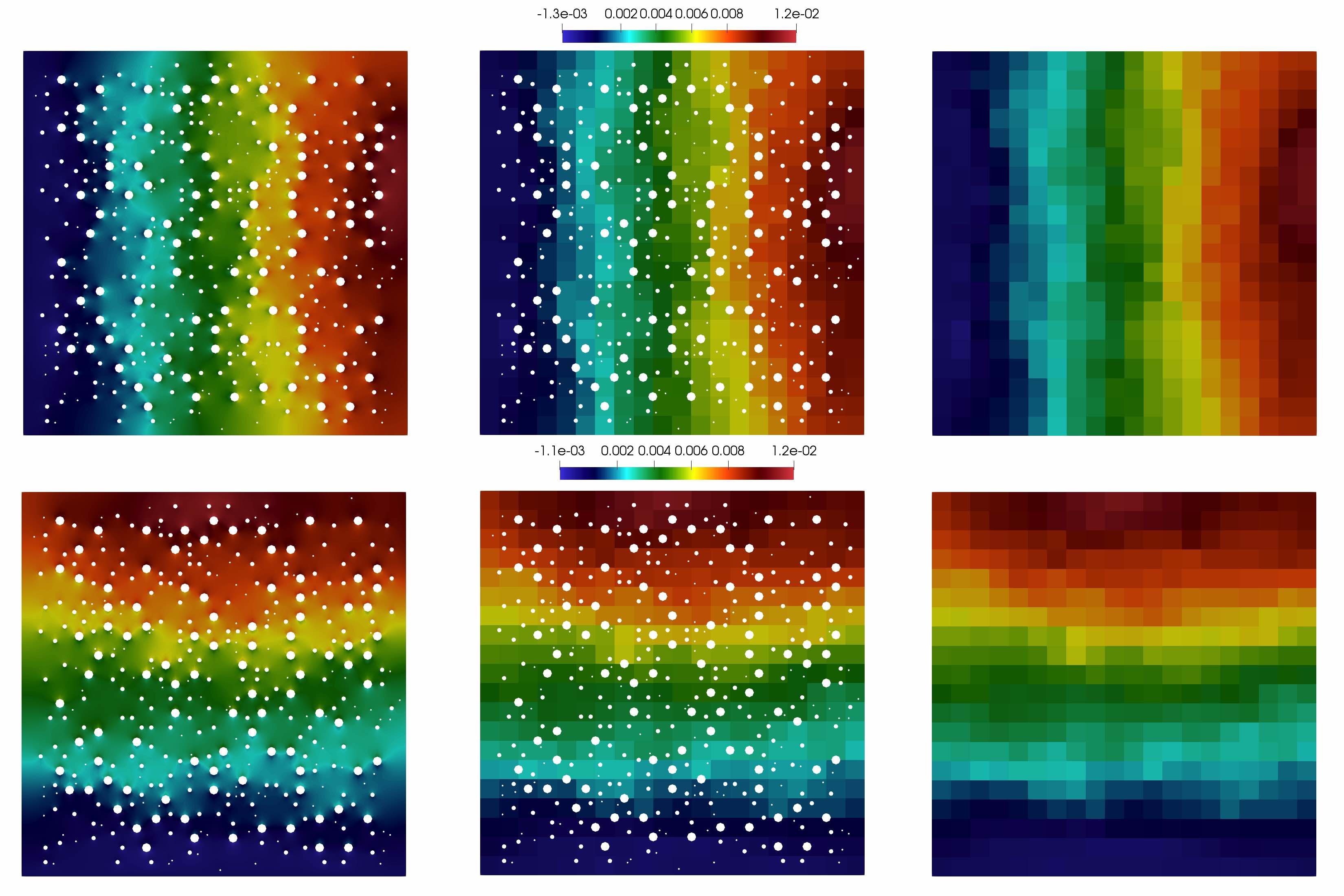}
\caption{Elasticity problem in the perforated domain shown in Figure~\ref{fig:mesh}. First column: fine-scale solution using $DOF_f = 30778$. Second column: mean value of fine-scale solution on a $20 \times 20$ coarse mesh. Third column: multiscale solutions computed on a $20 \times 20$ coarse mesh with $DOF_c = 1600$ using $4$ oversampling layers in basis constructions.
(The first row gives the first component of displacement and the second row gives the second component of displacement.)}
\label{fig:u20xy}
\end{figure}

\begin{table}[h!]
\centering
\begin{tabular}{ |c | c | c | }
\hline
$K^s$ & 
Type 1 & 
Type 2 \\
\hline
\multicolumn{3}{|c|}{Coarse mesh  $20 \times 20$} \\
\hline
$s = 1$ & 98.854	&	98.004 \\ \hline
$s = 2$ & 96.831	&	69.208 \\ \hline
$s = 3$ & 96.554	&	9.864 \\ \hline
$s = 4$ & 1.836		&	1.287 \\ \hline
\multicolumn{3}{|c|}{Coarse mesh  $40 \times 40$} \\
\hline
$s = 1$ & 99.792	&	99.820 \\ \hline
$s = 2$ & 97.716	&	97.768 \\ \hline
$s = 3$ & 91.475	&	79.359 \\ \hline
$s = 4$ & 24.594	&	24.329 \\ \hline
$s = 6$ & 0.637		&	0.642 \\ \hline
\end{tabular}
\,\,\,\,
\begin{tabular}{ |c | c | c | }
\hline
$K^s$ & 
$u_x$ & 
$u_y$ \\
\hline
\multicolumn{3}{|c|}{Coarse mesh  $20 \times 20$} \\
\hline
$s = 1$ & 95.451	&	96.073 \\ \hline
$s = 2$ & 77.983	&	73.635 \\ \hline
$s = 3$ & 10.026	&	13.585 \\ \hline
$s = 4$ & 1.959		&	0.928 \\ \hline
\multicolumn{3}{|c|}{Coarse mesh  $40 \times 40$} \\
\hline
$s = 1$ & 99.057	&	99.064 \\ \hline
$s = 2$ & 96.950	&	97.102 \\ \hline
$s = 3$ & 67.089	&	67.695 \\ \hline
$s = 4$ & 20.924	&	22.024 \\ \hline
$s = 6$ & 0.460		&	0.475 \\ \hline
\end{tabular}
\caption{Performance of our upscale method for the 
Laplace problem (left) and the elasticity problem (Right) in a perforated domain. Relative errors in percentage for the average solution on $20 \times 20$ and $40 \times 40$ coarse grids .}
\label{err-s}
\end{table}

For the model problems, we use following parameters:
\begin{itemize}
\item \textit{Laplace problem.} $k = 1$ and $f = 0$ with boundary conditions $u = 0$ for $x = 0$ and $y = 0$, $-k \nabla u \cdot n = 0$ for $x = 1$ and $y = 1$.
\item \textit{Elasticity problem}. $\mu = \frac{E}{2 (1 + \nu)}$, $\lambda = \frac{E \nu}{(1 + \nu) (1 - 2 \nu)}$, $E = 1$, $\nu = 0.3$, and $f = 0$ with boundary conditions $u_x = \sigma_y = 0$ for $x = 0$, $u_y = \sigma_x = 0$ for $y = 0$, $\sigma \, n = 0$ for $x = 1$ and $y = 1$.
\end{itemize}
On boundary of perforations, we set $-k \nabla u \cdot n = 1$ for Laplace problem and $\sigma \,  n = 1$ for elasticity problem.

Fine scale and upscaled solution are presented in Figure \ref{fig:u20} for Laplace problem and in Figure \ref{fig:u20xy} for elasticity problem.
In the first column, we give a fine-scale solution with $DOF_f = 15389$ for Laplace problem and $DOF_f = 30778$ for elasticity problem. In the second column, we show an average value of the fine-scale solution on a $20\times 20$ coarse grid. In the third column, we present a multiscale solutions computed on a $20 \times 20$ coarse grid with $4$ oversampling layers in the construction of basis functions. For the Laplace problem $DOF_c = 961$ and for the elasticity problem $DOF_c = 1600$. For both cases, the relative error for the solution corresponding to the background medium is about one percent.

In Table \ref{err-s}, we present relative errors 
%\marginpar{what norm?}
for both problems for two coarse grids and for different number of oversampling layers $K^s$ with $s = 1,2,3,4$ and $6$. From the numerical results, we observe a good convergence behaviour, when we take sufficient number of oversampled layers. For the coarse mesh with 400 cells, when we take 4 oversampling layers, we have $1.835 \%$ relative error for Type 1 basis functions and similar error for Type 2 basis functions for the Laplace problem. For the coarse mesh with 1600 cells, relative error is $0.637 \%$ for Type 1 basis and similar error for Type 2 basis. We note that, on the $20 \times 20$ coarse mesh, the size of upscaled system is $DOF_c = 731$ for Type 1 basis functions and $DOF_c = 961$ for Type 2 basis functions. For the $40 \times 40$ coarse mesh, we have $DOF_c = 2326$ for Type 1 basis functions and $DOF_c = 2354$ for Type 2 basis functions. We recall that the difference between Type 1 and Type 2 basis functions is the number of basis functions for perforations. For Type 1, when we have several perforations in coarse cell $K_i$, we use only one basis for handling all possible boundary conditions on perforations, but for Type 2, we use a basis for each perforation in a coarse cell.
For the elasticity problem, we present results for Type 2 basis functions and present errors for each component of solution $u_x$ and $u_y$, where for the coarse mesh with 1600 cells, we have less than one percent errors with 6 layers of oversampling. On the $20 \times 20$ coarse mesh, the size of upscaled system is $DOF_c = 1600$ and $DOF_c = 4000$ for the $40 \times 40$ coarse mesh. All results show good accuracy of the proposed method for both problems.

\subsection{Numerical results for time-dependent problem}\label{sec:num2}
% neumann
\begin{figure}[h!]
\centering
\includegraphics[width=1 \textwidth]{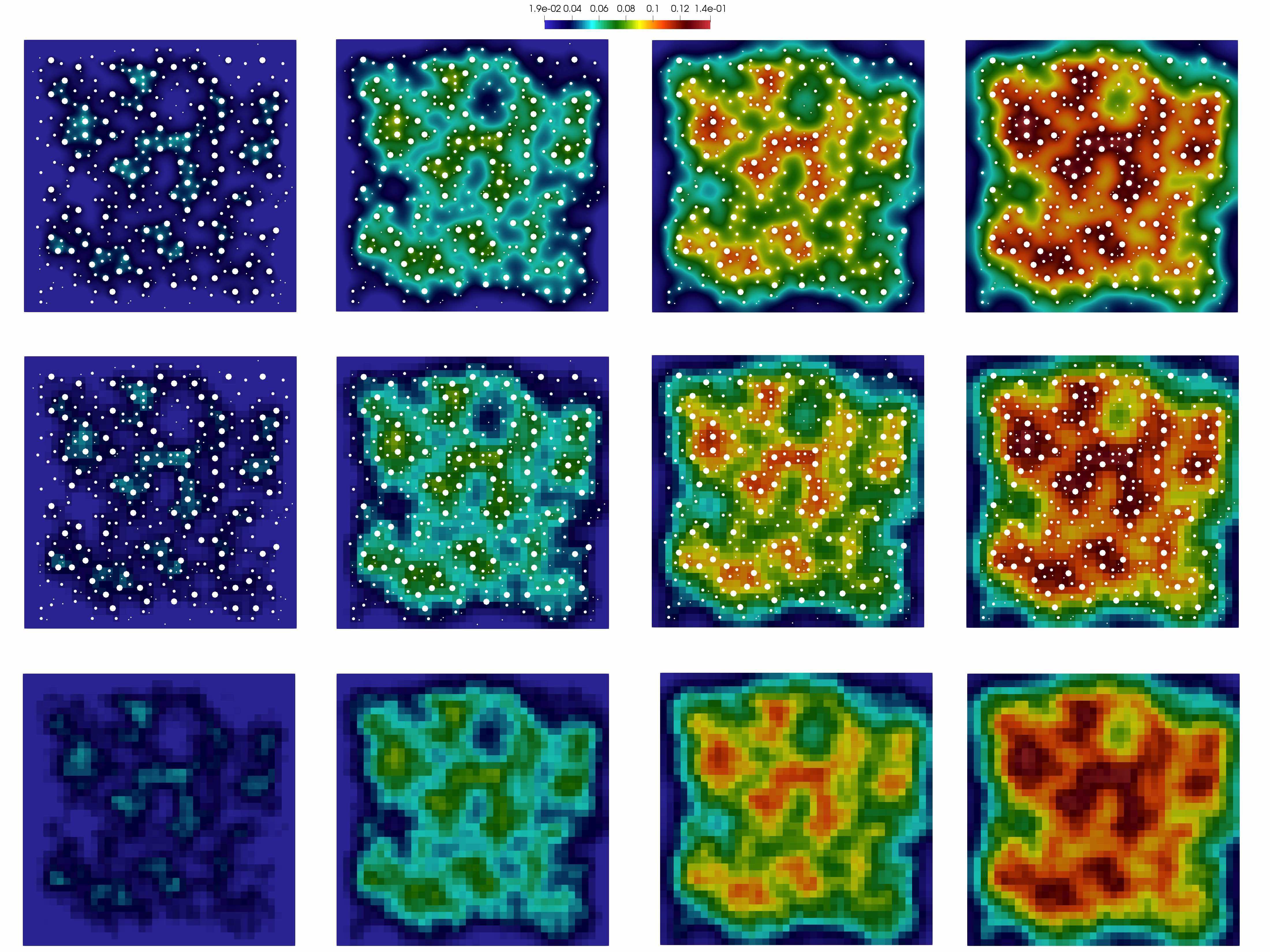}
\caption{Parabolic problem in a perforated domain with non-homogeneous Neumann boundary conditions on perforations (Test 1). First row: fine-scale solution with $DOF_f = 15389$. Second row: mean value of fine-scale solution on $40 \times 40$ coarse mesh. Third row: multiscale solutions $40 \times 40$ coarse mesh with $DOF_c = 2326$ and oversampling region with 6 coarse grid layers. We present the solution at different time instants $t_{5}   = 0.00125$, $t_{10} = 0.0025$, $t_{15} = 0.00375$  and $t_{20} = 0.005$ (from left to right).}
\label{fig:u-n}
\end{figure}

\begin{table}[h!]
\centering
\begin{tabular}{ |c | c | c | c | c | }
\hline
$K^s$ & 
$t_{5}   = 0.00125$ & 
$t_{10} = 0.0025$ & 
$t_{15} = 0.00375$ & 
$t_{20} = 0.005$ \\
\hline
\multicolumn{5}{|c|}{Coarse mesh  $20 \times 20$} \\
\hline
$s = 1$ & 3.865	&	3.581	&	3.468	&	3.399 \\ \hline
$s = 2$ & 3.429	&	3.324	&	3.302	&	3.261 \\ \hline
$s = 3$ & 3.412	&	3.318	&	3.308	&	3.278 \\ \hline
$s = 4$ & 2.735	&	1.553	&	1.061	&	0.798 \\ \hline
\multicolumn{5}{|c|}{Coarse mesh  $40 \times 40$} \\
\hline
$s = 1$ & 18.688	&	27.751 	&	38.829	&	47.005 \\ \hline
$s = 2$ & 1.570	&	1.433	&	1.389	&	1.390 \\ \hline
$s = 3$ & 1.361	&	1.265	&	1.194	&	1.129 \\ \hline
$s = 4$ & 0.866	&	0.453	&	0.308	&	0.239 \\ \hline
$s = 6$ & 0.862	&	0.443	&	0.304	&	0.224 \\ \hline
\end{tabular}
\caption{Parabolic problem in a perforated domain with non-homogeneous Neumann boundary conditions on perforations. Relative errors for the average value of solution on coarse grids $20 \times 20$ and $40 \times 40$. Type 1 basis functions are used in the simulations.}
\label{err-n1}
\end{table}

%\begin{table}[h!]
%\centering
%\begin{tabular}{ |c | c | c | c | c | }
%\hline
%$K^s$ & 
%$t_{5}   = 0.00125$ & 
%$t_{10} = 0.0025$ & 
%$t_{15} = 0.00375$ & 
%$t_{20} = 0.005$ \\
%\hline
%\multicolumn{5}{|c|}{Coarse mesh  $20 \times 20$} \\
%\hline
%$s = 1$ & 4.064	&	3.555	&	3.519	&	3.605 \\ \hline
%$s = 2$ & 2.678	&	1.531	&	1.071	&	0.837 \\ \hline
%$s = 3$ & 2.649	&	1.526	&	1.070	&	0.809 \\ \hline
%$s = 4$ & 2.634	&	1.522	&	1.060	&	0.794 \\ \hline
%\multicolumn{5}{|c|}{Coarse mesh  $40 \times 40$} \\
%\hline
%$s = 1$ & 34.963	&	67.030 	&	125.919	&	253.455 \\ \hline
%$s = 2$ & 2.458	&	2.745	&	2.879	&	2.998 \\ \hline
%$s = 3$ & 1.228	&	0.864	&	0.754	&	0.708 \\ \hline
%$s = 4$ & 0.866	&	0.453	&	0.308	&	0.238 \\ \hline
%$s = 6$ & 0.892	&	0.463	&	0.304	&	0.224 \\ \hline
%\end{tabular}
%\caption{Parabolic problem in perforated domain with non-homogeneous Neuman boundary conditions on perforations. Relative errors of the mean solution on a coarse grids $20 \times 20$ and $40 \times 40$. Type 2 basis functions}
%\label{err-n2}
%\end{table}

%Robyn boundary conditions

\begin{figure}[h!]
\centering
\includegraphics[width=1 \textwidth]{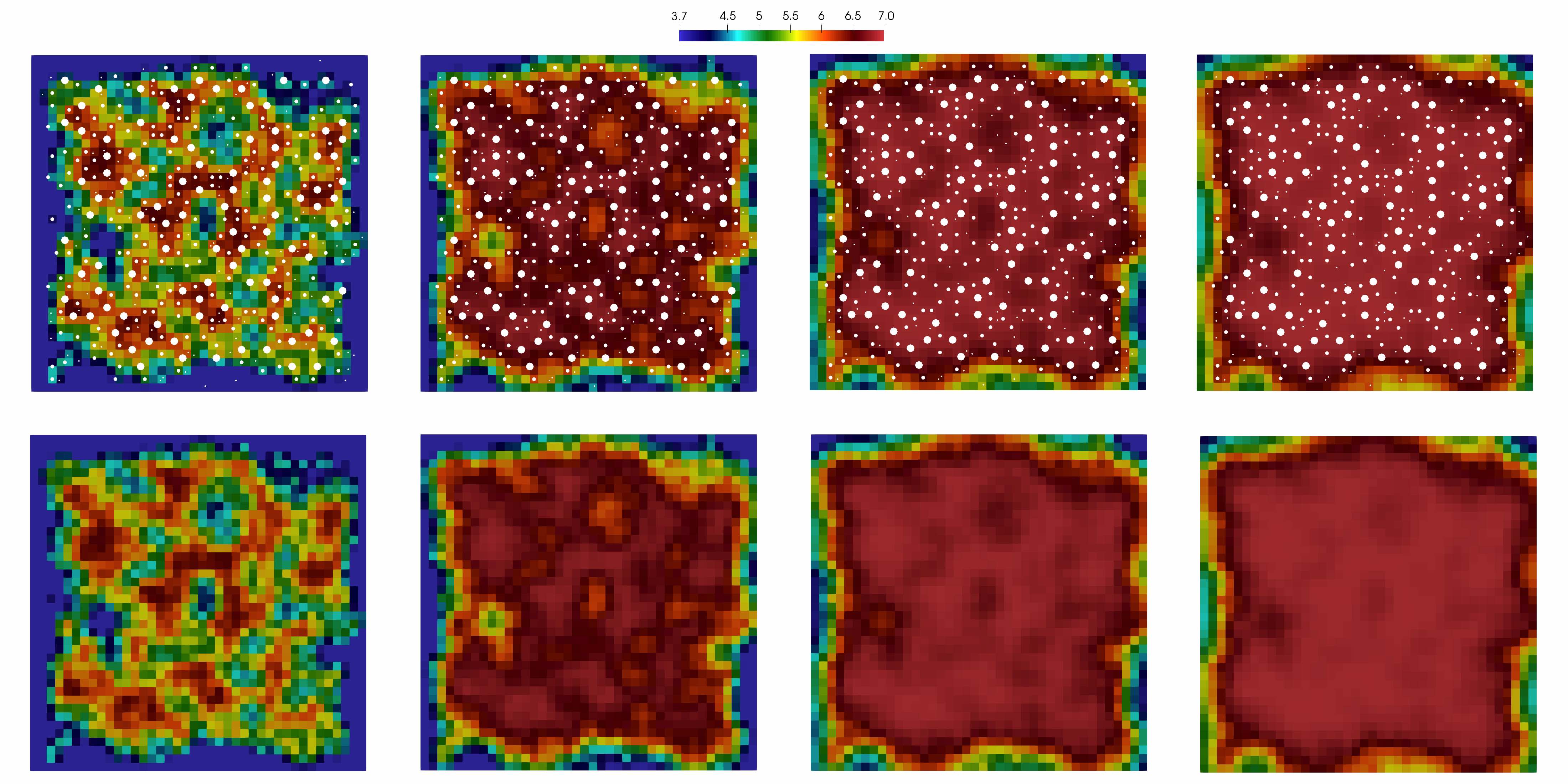}
\caption{Parabolic problem in a perforated domain with non-homogeneous Robin boundary conditions on perforations (Test 2). First row: fine-scale solution with $DOF_f = 15389$. Second row: mean value of fine-scale solution on $40 \times 40$ coarse mesh. Third row: multiscale solutions $40 \times 40$ coarse mesh with $DOF_c = 2326$ and oversampling region with 6 coarse grid layers. We present the solution at different time instants $t_{5}   = 0.00125$, $t_{10} = 0.0025$, $t_{15} = 0.00375$  and $t_{20} = 0.005$ (from left to right).}
\label{fig:u-r}
\end{figure}

\begin{figure}[h!]
\centering
\includegraphics[width=1 \textwidth]{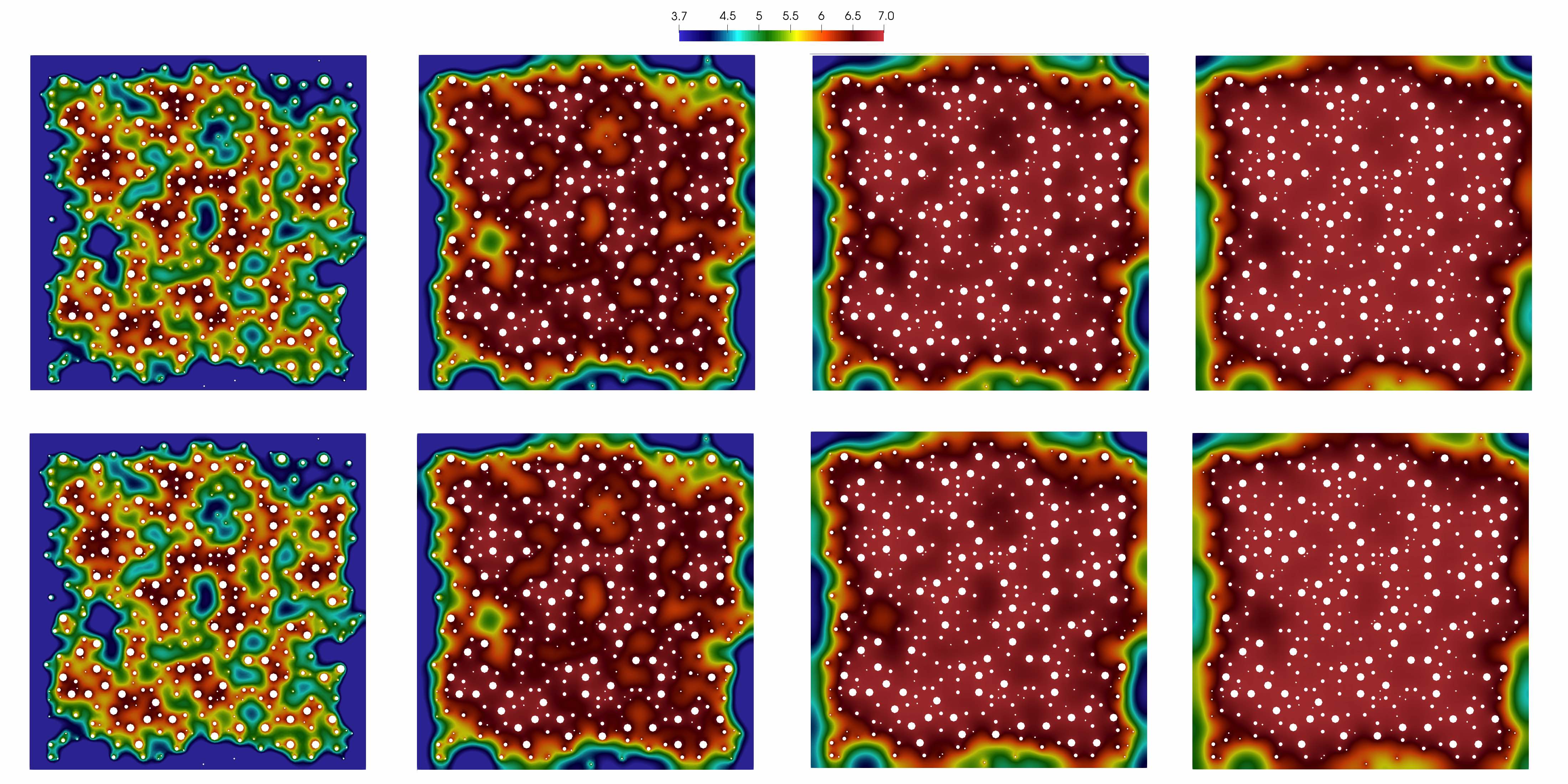}
\caption{Parabolic problem in a perforated domain with non-homogeneous Robin boundary conditions on perforations. First row: fine-scale solution with $DOF_f = 15389$. Second row: downscale solutions the $40 \times 40$ coarse mesh with $DOF_c = 2326$ using 6 oversampling layers in basis construction. We present the solution at different time instants $t_{5}   = 0.00125$, $t_{10} = 0.0025$, $t_{15} = 0.00375$  and $t_{20} = 0.005$ (from left to right).}
\label{fig:u-r2}
\end{figure}

\begin{table}[h!]
\centering
\begin{tabular}{ |c | c | c | c | c | }
\hline
$K^s$ & 
$t_{5}   = 0.00125$ & 
$t_{10} = 0.0025$ & 
$t_{15} = 0.00375$ & 
$t_{20} = 0.005$ \\
\hline
\multicolumn{5}{|c|}{Coarse mesh  $20 \times 20$} \\
\hline
$s = 1$ & 21.897	&	25.518	&	27.812	&	29.375 \\ \hline
$s = 2$ & 16.030	&	17.217	&	17.928	&	18.397 \\ \hline
$s = 3$ & 15.844	&	16.938	&	17.588	&	17.991 \\ \hline
$s = 4$ & 1.948		&	1.199	&	0.938	&	0.806 \\ \hline
\multicolumn{5}{|c|}{Coarse mesh  $40 \times 40$} \\
\hline
$s = 1$ & 50.837	&	60.026 	&	64.118	&	66.390 \\ \hline
$s = 2$ & 11.800	&	13.691	&	15.129	&	16.210 \\ \hline
$s = 3$ & 8.818	&	8.632	&	8.550	&	8.494 \\ \hline
$s = 4$ & 0.758	&	0.449	&	0.335	&	0.280 \\ \hline
$s = 6$ & 0.738	&	0.442	&	0.332	&	0.277 \\ \hline
\end{tabular}
\caption{Parabolic problem in a perforated domain with non-homogeneous Neumann boundary conditions on perforations (Test 1). Relative errors for the average of the solution on coarse grids with sizes $20 \times 20$ and $40 \times 40$. Type 1 basis functions are used in the simulations.}
\label{err-r1}
\end{table}

\begin{table}[h!]
\centering
\begin{tabular}{ |c | c | c | c | c | }
\hline
$K^s$ & 
$t_{5}   = 0.00125$ & 
$t_{10} = 0.0025$ & 
$t_{15} = 0.00375$ & 
$t_{20} = 0.005$ \\
\hline
\multicolumn{5}{|c|}{Coarse mesh  $20 \times 20$} \\
\hline
$s = 1$ & 12.609	&	15.717	&	17.912	&	19.466 \\ \hline
$s = 2$ & 2.253	&	1.470	&	1.236	&	0.166 \\ \hline
$s = 3$ & 2.067	&	1.241	&	0.932	&	0.771 \\ \hline
$s = 4$ & 2.059	&	1.237	&	0.931	&	0.770 \\ \hline
\multicolumn{5}{|c|}{Coarse mesh  $40 \times 40$} \\
\hline
$s = 1$ & 54.461	&	63.396 	&	67.208	&	69.288 \\ \hline
$s = 2$ & 8.568	&	11.829	&	13.889	&	15.278 \\ \hline
$s = 3$ & 1.289	&	1.327	&	1.484	&	1.622 \\ \hline
$s = 4$ & 0.760	&	0.450	&	0.336	&	0.280 \\ \hline
$s = 6$ & 0.740	&	0.440	&	0.331	&	0.274 \\ \hline
\end{tabular}
\caption{Parabolic problem in a perforated domain with non-homogeneous Neumann boundary conditions on perforations. Relative errors for the average of the solution on coarse grids with sizes $20 \times 20$ and $40 \times 40$. Type 2 basis functions are used in the simulations.}
\label{err-r2}
\end{table}

Next, we consider the time-dependent problem. We perform numerical simulations on the same perforated domain depicted in Figure~\ref{fig:mesh} and use similar coarse grids.
In addition, we use the following parameters: $c = 1$, $k = 1$ and $f = 0$ with boundary conditions $-k \nabla u \cdot n = 0$ on $\partial \Omega \backslash \Gamma$. We consider two test cases, where we set $- k \nabla u \cdot n = 1$ (test 1) and $- k \nabla u \cdot n = 100 (u - 7)$ (test 2) on boundary of perforations. We consider $T_{max} = 0.005$ and use 20 time steps.

In Figure \ref{fig:u-n}, we present the fine scale and upscaled solutions for non-homogeneous Neumann boundary conditions (test 1) and in Figure \ref{fig:u-r} for non-homogeneous Robin boundary conditions (test 2). The downscale solution is shown in Figure \ref{fig:u-r2} for test 2. The size of the fine grid system is $DOF_f = 15389$. Coarse scale system has size $DOF_c = 731$ for Type 1 basis functions and $DOF_c = 961$ for Type 2 basis functions on the $20 \times 20$ coarse grid. For the $40 \times 40$ coarse mesh, we have $DOF_c = 2326$ for Type 1 basis functions and $DOF_c = 2354$ for Type 2 basis functions. From the Figures~\ref{fig:u-n}-\ref{fig:u-r2}, we observe very good agreement between the fine-scale solution and the computed upscaled solution. 

In Table \ref{err-n1}, we present relative errors for two choices of coarse grids and for different number of oversampling layers $K^s$ with $s = 1,2,3,4$ and $6$ for non-homogeneous Neumann boundary conditions (test 1). We used Type 1 multiscale basis functions and observe a good convergence behaviour, when we take a sufficient number of oversampled layers. For example, for coarse mesh with 1600 cells, when we take 4 oversampling layers, we have less the one percent relative error. Finally,
in Tables \ref{err-r1} and \ref{err-r2}, we present relative errors for Type 1 and 2 multiscale basis functions. The results show good accuracy of the proposed method for Type 1 and 2 basis functions, but for Type 2 we can take smaller number of oversampling layers. For the coarse mesh with 1600 cells, we have less than one percent relative errors with 4 layers of oversampling.

\section{Conclusion}\label{sec:con}

We presented an upscaling method for problems in perforated domains with non-homogeneous boundary conditions on perforations. In this method, we construct multiscale basis function for background medium and additional multiscale basis for perforations, that help to handle non-homogeneous boundary conditions.
We proposed a method and presented numerical results for Laplace, elasticity and parabolic problems. Numerical results show that the proposed method can provide good accuracy and give a significant reduction of the size of system for problems in perforated domains.
The resulting upscaled model has minimal size and the computed solution has a physical meaning on the coarse grid.

\section{Acknowledgements}
MV's  and DS's works are supported by the mega-grant of the Russian Federation Government (N 14.Y26.31.0013). 
EC's work is partially supported by Hong Kong RGC General Research Fund (Project 14304217)
and CUHK Direct Grant for Research 2016-17.

\bibliographystyle{plain}
\bibliography{lit}

\end{document}